\documentclass[psamsfonts,fceqn,leqno]{amsart}
\usepackage{amsthm}
\usepackage{amsmath}
\usepackage{amssymb}
\usepackage{color}
\usepackage{mathrsfs}

\newtheorem{theorem}{Theorem}[section]
\newtheorem{corollary}[theorem]{Corollary}

\newtheorem{lemma}[theorem]{Lemma}
\newtheorem{question}[theorem]{Question}
\newtheorem{problem}[theorem]{Problem}
\theoremstyle{definition}

\newtheorem{remark}[theorem]{Remark}

\def\N{{\mathbb N}}

\numberwithin{equation}{section}

  \title[Pairwise $k$-Semi-Stratifiable Bispaces and Topological Ordered Spaces]
  {Pairwise $k$-Semi-Stratifiable Bispaces and Topological Ordered Spaces}

  \author{Kedian Li}
  \address{\emph{Kedian Li}: School of Mathematics and Statistics, Minnan Normal University, Zhangzhou 363000, 
  P. R. China}
  \email{likd56@126.com}
  \author{Jiling Cao$^\ast$}
  \address{\emph{Jiling Cao}: School of Engineering, Computer and Mathematical Sciences, Auckland
  University of Technology, Private Bag 92006, Auckland 1142, New Zealand}
  \email{jiling.cao@aut.ac.nz}

  \thanks{The first-author is supported by the National Nature Science Fundation of China, grant 
  No. 61379021,11471153) and the Natural Science Foundation of Fujian Province, China, grant No. 2013J01029).
  The second author thanks the support of the National Natural Science Foundation of China, grant No. 
  11571158, and the paper was written when he visited Minnan Normal University in April 2016 as 
  Min Jiang Scholar Guest Professor.}
  
  \thanks{*Corresponding author}


  
\begin{document}
  	
  \begin{abstract}
  In this paper, we continue to study pairwise ($k$-semi-)stratifiable bitopological spaces. Some new 
  characterizations of pairwise $k$-semi-stratifiable bitopological spaces are provided. Relationships 
  between pairwise stratifiable and pairwise $k$-semi-stratifiable bitopological spaces are 
  further investigated, and an open question recently posed by Li and Lin in \cite{LL} is completely 
  solved. We also study the quasi-pseudo-metrizability of a topological ordered space $(X, \tau, 
  \preccurlyeq)$. It is shown that if $(X, \tau, \preccurlyeq)$ is a ball transitive topological 
  ordered $C$- and $I$-space such that $\tau$ is metrizable, then its associated bitopological space 
  $(X,\tau^{\flat},\tau^{\natural})$ is quasi-pseudo-metrizable. This result provides a partial
  affirmative answer to a problem in \cite{KM}. 
  \end{abstract}

  \maketitle
  
  \section{Introduction}
  
  Undoubtedly, topology and order are not only important topics in mathematics but also applicable in 
  many other disciplines. For example, nonsymmetric notions of distance are needed for mathematical 
  modelling in the natural, physical and cybernetic sciences and the corresponding topological notion is
  that of a quasi-metric or a quasi-pseudo-metric. The study of quasi-metrizable spaces naturally leads
  to the concepts of quasi-uniformities and bitopological spaces. In this aspect, Kelly's seminal paper
  \cite{Ke} made pioneer contributions. On the other hand, the notions of a sober space and the Scott
  topology, in align with the investigation of partial orders and pre-orders, are useful in theoretical
  computer science in the study of algorithms which act on other algorithms. Moreover, finite topological 
  spaces (i.e., finite pre-orders) can be used to construct a mathematical model of a video monitor screen 
  which may be useful in computer graphics.
    
  \medskip
  Since Kelly's work in \cite{Ke}, bitopological spaces have attracted the attention of many researchers.
  For example, Reilly \cite{Rei:72} explored separation axioms for bitopological spaces, Cooke and Reilly
  \cite{CR:75} discussed the relationships between six definitions of bitopological compactness appeared 
  in the literature, Raghavan and Reilly \cite{RaRe:86} introduced a notion of bitopological 
  paracompactness and established a bitopological version of Michael's classical characterization of 
  regular paracompact spaces. In addition to separation and covering properties, generalized metric 
  properties have also been considered in the setting of bitopological spaces. In this direction, Fox
  \cite{Fox} discussed the quasi-metrizability of bitopological spaces, pairwise stratifiable bitopological
  spaces and their generalizations have been introduced and studied by \cite{GR}, \cite{MR:91} and 
  \cite{MMR:86}.
  
  \medskip
  Interplay between topology and order has been a very interesting area. In 1965, Nachbin's book \cite{Na} 
  was published. This book is one of general references on the subject available today, and it covers results 
  obtained by the author in his research on spaces with structures of order and topology. Among many topics
  in this line, McCartan \cite{Mc} studied bicontinuous (herein called $C$-space and $I$-space) pre-ordered 
  topological spaces and investigated the relationships between the topology of such a space and two 
  associated convex topologies, and Faber \cite{Fa} studied metrizability in generalized ordered spaces. 
  Recently, there have been some renewed interests in the study of generalized metric properties in 
  bitopological spaces and topological ordered spaces. K\"{u}nzi and Mushaandja investigated the
  quasi-pseudo-metrizability of a topological ordered space in \cite{KM} and \cite{Mu}, respectively. 
  They obtained some results related to the upper topology $\tau^{\natural}$ and the lower topology 
  $\tau^{\flat}$ of a metrizable ordered space $(X, \tau, \preccurlyeq)$ which is both a $C$- and an 
  $I$-space in the sense of Priestle \cite{PR}. Moreover, Li \cite{Li:10} as well as Li and Lin 
  \cite{Kedian-Fucai:13} carried on the study of pairwise (semi-)stratifiable bispaces and established 
  some new characterizations for these classes of bitopological spaces. In a very recent paper \cite{LL}, 
  Li and Lin further introduced and studied the class of $k$-semi-stratifiable bitopological spaces.
  
  \medskip
  The main purpose of this paper is to continue the study of pairwise $k$-semi-stratifiable bitopological 
  spaces and their relationships with and applications to the quasi-pseudo-metrizability of a topological 
  ordered spaces. For the sake of self-completeness, in Section \ref{sec:notation}, we introduce the 
  necessary definitions and terminologies. In Section \ref{sec:charac}, we provide some new 
  characterizations of pairwise $k$-semi-stratifiable bitopological spaces in terms of $g$-functions， 
  $\sigma$-cushioned pair $k$-networks and $cs$-networks. In Section \ref{sec:relation}, we consider 
  some conditions under which a pairwise $k$-semi-stratifiable bitopological space is pariwise 
  stratifiable. An open question posed in \cite{LL} is completely solved and results in \cite{Cao:99} 
  and \cite{Cao:05} are extended to the setting of bitopological spaces. In the last section, we consider 
  the quasi-pseudo-metrizability of a topological ordered spaces and provide a partial affirmative answer 
  to an open problem of K\"{u}nzi and Mushaandja in \cite{KM}. 
  
  \medskip
  Our notations in this paper are standard. For any undefined concepts and terminologies, we refer the
  reader to \cite{E} or \cite{Gr}.

  \section{Preliminaries and notations} \label{sec:notation}
  
  A \emph{quasi-pseudo-metric} $d$ on a nonempty set $X$ is a non-negative real-valued function $d: X 
  \times X \to \mathbb R_+$ such that (i) $d(x,x)=0$ and (ii) $d(x,z) \le d(x,y) + d(y,z)$, for all 
  $x, y, z \in X$. If $d$ is a quasi-pseudo-metric on $X$, then the ordered pair $(X,d)$ is called a 
  \emph{quasi-pseudo-metric space}. Every quasi-pseudo-metric $d$ on $X$ induces a topology 
  $\tau(d)$ on $X$ which has as a base the family $\{B_d(x, \epsilon): x\in X, \epsilon >0\}$, where 
  $B_d(x,\epsilon) =\{ y \in X: d(x,y) < \epsilon\}$. Every quasi-pseudo-metric $d$ on $X$ induces a 
  conjugate quasi-pseudo-metric $d^{-1}$ on $X$, defined by $d^{-1}(x,y) =d(y,x)$ for all $x, y \in X$. 
  A \emph{bitopological space} \cite{Ke} (for short, \emph{bispace} \cite{GR}) is a triple $(X, \tau_1, 
  \tau_2)$, where $X$ is a nonempty set, topologies $\tau_1$ and $\tau_2$ are two topologies on $X$. A 
  bispace $(X, \tau_1, \tau_2)$ is called \emph{quasi-pseudo-metrizable}, if there is a quasi-pseudo-metric 
  $d$ on $X$ such that $\tau(d) =\tau_1$ and $\tau(d^{-1}) =\tau_2$.
  
  \medskip
  Let $(X, \tau_1, \tau_2)$ be a bispace. For $i=1,2$, let ${\mathscr F}_i(X)$ denote the family of all 
  $\tau_i$-closed subsets of $X$. For a subset $A$ of $X$, let $\overline{A}^{\tau_i}$ and ${\rm int}_{\tau_i}(A)$ 
  denote the closure and interior of $A$ with respect to $\tau_i$, $i=1,2$, respectively. For $i, j=1,2$ with 
  $i\ne j$, $(X, \tau_1, \tau_2)$ is called \emph{$\tau_i$-semi-stratifiable with respect to $\tau_j$} if there 
  exists an operator $G_{ij}: \mathbb{N} \times {\mathscr F}_i(X) \rightarrow \tau_j$ satisfying 
  (i) $H=\bigcap_{n\in \mathbb N} G_{ij}(n,H)$ for all $H\in {\mathscr F}_i(X)$, 
  (ii) if $H, K \in {\mathscr F}_i(X)$ with $H\subseteq K$, then $G_{ij}(n,H)\subseteq G_{ij}(n,K)$ for 
  all $n\in \mathbb N$. Furthermore, if $G_{ij}$ satisfies (ii) and (i)'  $H=\bigcap_{n\in \mathbb N} 
  \overline{G_{ij}(n,H)}^{\tau_i}$ for all $H\in {\mathscr F}_i(X)$, 
  then $(X, \tau_1, \tau_2)$ is called \emph{$\tau_i$-stratifiable with respect to $\tau_j$}. Moreover, 
  $(X, \tau_1, \tau_2)$ is called \emph{pairwise (semi-)stratifiable} \cite{Fox}, \cite{GR} and \cite{MMR:86}, 
  if it is both $\tau_1$-semi-stratifiable with respect to $\tau_2$ and $\tau_2$-semi-stratifiable with respect 
  to $\tau_1$.
  
  \medskip
  Recently, Li and Lin \cite{LL} introduced the concept of a pairwise $k$-semi-stratifiable bispace, which is a 
  natural extension of a $k$-semi-straitifiable space introduced in \cite{Lutzer:71} to the setting of bispaces. 
  For $i, j=1,2$ with $i\ne j$, a bispace $(X, \tau_1, \tau_2)$ is called \emph{$\tau_i$-$k$-semi-stratifiable 
  with respect to $\tau_j$} if there exists an operator $G_{ij}: \mathbb{N} \times {\mathscr F}_i(X) \rightarrow 
  \tau_j$ satisfying 
  (i) $H=\bigcap_{n\in \mathbb N} G_{ij}(n,H)$ for all $H\in {\mathscr F}_i(X)$, 
  (ii) if $H,K\in {\mathscr F}_i(X)$ with $H\subseteq K$, then $G_{ij}(n,H)\subseteq G_{ij}(n,K)$ for 
  all $n\in \mathbb N$, 
  (iii) if $K\subseteq X$ is $\tau_i$-compact and $H\in {\mathscr F}_i(X)$ such that $H \cap K =
  \emptyset$, then $K\cap G_{ij}(n,H)=\emptyset$ for some $n\in\mathbb{N}$. 
  In addition, $(X, \tau_1, \tau_2)$ is called \emph{pairwise $k$-semi-stratifiable} \cite{LL} if it is 
  both $\tau_1$-$k$-semi-stratifiable with respect to $\tau_2$ and $\tau_2$-$k$-semi-stratifiable with respect 
  to $\tau_1$. 
  
  \medskip
  The next lemma, which gives an important dual characterization of pairwise $k$-semi-stratifiable bispaces, 
  will be used in the sequel.
  
  \begin{lemma}[\cite{LL}] \label{lem:pksemi}
  A bispace $(X, \tau_1, \tau_2)$ is pairwise $k$-semi-stratifiable if, and only if, for any $i, j=1,2$ with 
  $i\ne j$, there is an operator $F_{ij}: \mathbb{N}\times \tau_i \rightarrow {\mathscr F}_j(X)$ satisfying
  \begin{enumerate}
  \item[{\rm (1)}] $U=\bigcup_{n\in \mathbb{N}}F_{ij}(n,U)$ for all $U \in \tau_i$;
  		
  \item[{\rm (2)}] if $U, V\in \tau_i$ with $U \subseteq V$, then $F_{ij}(n,U)\subseteq F_{ij}(n,V)$ for all 
  $n\in \mathbb{N}$;
  		
  \item[{\rm (3)}] if $K \subseteq X$ is $\tau_i$-compact and $U\in \tau_i$ with $K \subseteq U$, then 
  $K\subseteq F_{ij}(n,U)$ for some $n\in\mathbb{N}$.
  \end{enumerate}
  In addition, the operator $F_{ij}$ can be required to be monotone with respect to $n$, that is, $F_{ij}(n,U) 
  \subseteq F_{ij}(n+1,U)$ for all $n\in \mathbb{N}$ and all $U\in \tau_i$. 
  \end{lemma}
  
  By definition, every pairwise stratifiable bispace is pairwise $k$-semi-straitifable, and every pairwise 
  $k$-semi-stratifiable bispace is pairwise semi-stratifiable. Recall that a bispace $(X, \tau_1, \tau_2)$ 
  is said to be \emph{pairwise monotonically normal} \cite{MR:91} if to each pair $(H,K)$ of disjoint subsets 
  of $X$ such that $H\in {\mathscr F}_i(X)$ and $K\in {\mathscr F}_j(X)$ ($i,j=1,2$ and $i \ne j$), we can 
  assign a set $D_{ij}(H,K) \in \tau_j$ such that 
  (i) 
  \[
  H \subseteq D_{ij}(H,K) \subseteq \overline{D_{ij}(H,K)}^{\tau_i} \subseteq X\smallsetminus K,
  \] 
  (ii) if the pairs $(H,K)$ and $(H',K')$ satisfy $H \subseteq H'$ and $K' \subseteq K$, then 
  $D_{ij}(H,K) \subseteq D_{ij}(H', K')$.  
  
  \medskip
  The following result, established by Mar\'{i}n and Romaguera in \cite{MR:91}, extends the celebrated 
  result of Heath et al. in \cite{HLZ:73} on monotonically normal spaces.
  
  \begin{theorem}[\cite{MR:91}] \label{thm:bidecom}
  A bispace $(X, \tau_1, \tau_2)$ is pairwise stratifiable if, and only if, it is a pairwise monotonically 
  normal and pairwise semi-stratifiable bispace.
  \end{theorem}
  
  \begin{corollary} \label{coro:bidecom}
  A pairwise monotonically normal and pairwise $k$ semi-stratifiable bispace is pairwise stratifiable.
  \end{corollary}
  
  A \emph{topological ordered space} $(X, \tau, \preccurlyeq)$ is a nonempty set $X$ endowed with a topology 
  $\tau$ and a partial order $\preccurlyeq$. A subset $A$ of $X$ is said to be an \emph{upper set} of $X$ 
  if $x \preccurlyeq y$ and $x \in  A$ imply that $y \in  A$. Similarly, we say that a subset $A$ of $X$ is 
  a \emph{lower set} of $X$ if $y \preccurlyeq x$ and $x \in A$ imply that $y \in A$. Let $\tau^{\flat}$ 
  denote the collection of $\tau$-open lower sets of $X$ and $\tau^{\natural}$ denote the collection of 
  $\tau$-open upper sets of $X$. Then, $\tau^{\flat}$ and $\tau^{\natural}$ are two topologies on $X$ and 
  thus $(X,\tau^{\flat}, \tau^{\natural})$ is a bispace. For any subset $A$ of $X$, $i(A)$ (resp. $d(A)$) 
  will denote the intersection of all upper (lower) sets of $X$ containing $A$. Note that $i(A)$ (resp. 
  $d(A))$ is the smallest upper (resp. lower) set containing $A$. It is easy to see that $A = i(A)$ if, 
  and only if, $A$ is an upper set. Similarly, $A = d(A)$ if, and only if, $A$ is a lower set. Following 
  Priestley \cite{PR}, we recall that a topological ordered space $(X,\tau,\preccurlyeq)$ is said to be 
  a \emph{$C$-space} if $d(F)$ and $i(F)$ are closed whenever $F$ is a closed subset of $X$. Similarly, a 
  topological ordered space $(X,\tau,\preccurlyeq)$ is called an \emph{$I$-space} if $d(G)$ and $i(G)$ 
  are open whenever $G$ is an open subset of $X$. 
  
  \section{Some new characterizations of pairwise\\ $k$-semi-stratifiable bispaces}\label{sec:charac}
  
  In \cite{LL}, Li and Lin characterized pairwise $k$-semi-stratifiable bispaces in terms of pairwise
  $g$-functions and extensions of semi-continuous functions. In this section, we continue to investigate
  how to characterize pairwise $k$-semi-stratifiable bispaces. Our first two results can be regarded as 
  either improvements or extensions of a theorem in \cite{LL}. In addition, we also use cushioned pair 
  $k$-networks and $cs$-networks to characterize pairwise $k$-semi-stratifiable bispaces. 
  
  \medskip
  Let $(X, \tau_1, \tau_2)$ be a bispace. A \emph{pairwise $g$-function} on $(X, \tau_1, \tau_2)$ is a 
  pair of functions $(g_1, g_2)$ such that for $i=1, 2$, $g_i: \mathbb N \times X \to \tau_i$ satisfies 
  $x\in g_i(n, x)$ and $g_i(n+1, x) \subseteq g_i(n,x)$ for all $n \in \mathbb N$ and $x\in X$. 
  A pairwise family ${\mathscr B} = \left\{\left(B_\alpha^1,B_\alpha^2 \right): \alpha \in \Delta \right
  \}$ of subsets of $X$ is called \emph{$\tau_j$-cushioned}, where $j=1,2$, if for any $\Delta' \subseteq 
  \Delta$, 
  \[
  \overline{\bigcup \left\{B_\alpha^1: \alpha \in \Delta' \right\}}^{\tau_j} \subseteq \bigcup \left\{B_\alpha^2: 
  \alpha \in \Delta' \right\}.
  \]
  Furthermore, if ${\mathscr B}$ is a countable union of $\tau_j$-cushioned families, then it is called 
  \emph{$\sigma$-$\tau_j$-cushioned}. A pairwise family ${\mathscr B} = \left\{ \left(B_\alpha^1, B_\alpha^2 
  \right): \alpha \in \Delta\right\}$ is called \emph{a pair $\tau_i$-$k$-network} if for any $\tau_i$-compact 
  set $K$ and any set $U\in \tau_i$ with $K \subseteq U$, there is a finite subset $\Delta'$ of $\Delta$ 
  such that
  \[
  K \subseteq \bigcup \{B_\alpha^1: \alpha \in \Delta'\} \subseteq \bigcup \{B_\alpha^2: \alpha \in \Delta\}
  \subseteq U.
  \]
  
  \medskip
  Our first result improves the equivalence of (1) and (2) in \cite[Theorem 2.1]{LL}.
  
  \begin{theorem} \label{thm:pksemig1}
  Let $(X, \tau_1, \tau_2)$ be a bispace such that $(X, \tau_i)$ is $T_1$-space for $i=1,2$. Then 
  $(X, \tau_1, \tau_2)$ is pairwise $k$-semi-straitifiable if, and only if, there is a pairwise $g$-function 
  $(g_1, g_2)$ such that for $i,j =1,2$ with $i\ne j$, if $K$ is a $\tau_i$-compact set and $H$ is a 
  $\tau_i$-closed set with $K \cap H =\emptyset$, then 
  \[
  K \cap \left(\bigcup \{ g_j(m, x): x\in H\} \right)= \emptyset
  \]
  for some $m \in \mathbb N$.
  \end{theorem}
  
  \begin{proof}
  Necessity. Suppose that $(X, \tau_1, \tau_2)$  is pairwise $k$-semi-stratifiable. For $i,j =1,2$ and 
  $i\ne j$, let $G_{ij}: \mathbb{N} \times {\mathscr F}_i(X) \rightarrow \tau_j$ be an operator satisfying
  the definition of a pairwise $k$-semi-stratifiable bispace. Without loss of generality, $G_{ij}$ 
  can be required to be monotone with respect to $n$. Define a function $g_j: \mathbb N \times X \to \tau_j$ 
  such that $g_j(n, x) = G_{ij}(n, \{x\})$ for all $n \in \mathbb N$ and $x\in X$. Clearly, $(g_1,g_2)$ is a 
  pairwise $g$-function. If $K$ is a $\tau_i$-compact subset and $H$ is a $\tau_i$-closed subset with $K \cap 
  H = \emptyset$, then $K \cap G_{ij}(m, H) = \emptyset$ for some $m \in \mathbb N$. Note that 
  \[
  \bigcup \{ g_j(m, x): x\in H\}\subseteq G_{ij}(m, H).
  \] 
  It follows that 
  \[
  K \cap \left(\bigcup \{ g_j(m, x): x\in H\} \right)= \emptyset.
  \]
	
  Sufficiency. Let $(g_1, g_2)$ be a pairwise $g$-function satisfying the assumption in the theorem.   
  For each $\tau_i$-closed subset $H$ and $n \in \mathbb N$, define
  \[
  G_{ij}(n, H) = \bigcup \{g_j(n, x): x \in H\}.
  \]
  We shall verify that $G_{ij}$ is an operator satisfying conditions (i) in the definition of a pairwise 
  $k$-semi-stratifiable bitopological space, as (ii) and (iii) hold trivially. It is clear that $H \subseteq 
  \bigcap_{n \in \mathbb N} G_{ij}(n, H)$. If $p \not \in H$, as $\{p \}$ is compact, then the assumption in
  the theorem implies that there must be some $m \in \mathbb N$ such that $p \not \in \bigcup \{g_j(m, x): x 
  \in H\}$. It follows that $p \not \in G_{ij}(m, H)$. Thus, $H = \bigcap_{n \in \mathbb N} G_{ij}(n, H)$.
  \end{proof}
  
  \begin{theorem} \label{thm:pksemig2}
  Let $(X, \tau_1, \tau_2)$ be a bispace such that $(X, \tau_i)$ is Hausdorff for $i=1,2$. 
  Then the following statements are equivalent.
  \begin{itemize}
  \item[{\rm (1)}] $(X, \tau_1, \tau_2)$ is pairwise $k$-semi-straitifiable.
		
  \item[{\rm (2)}] There is a pairwise $g$-function $(g_1, g_2)$ such that for $i,j =1,2$ with $i\ne j$, if 
  $\{x_n: n \in \mathbb N\}$ is a sequence $\tau_i$-convergent to $p$ and $H$ is a $\tau_i$-closed subset 
  with $\left(\{ p \} \cup \{x_n: n \in \mathbb N\}\right)\cap H = \emptyset$, then 
  \[
  \left(\{p \} \cup \{x_n: n \in \mathbb N\}\right) \cap\left(\bigcup \{ g_j(m, x): x\in H\}\right)= 
  \emptyset
  \] 
  for some $m \in \mathbb N$.
		
  \item[{\rm (3)}] There is a pairwise $g$-function $(g_1, g_2)$ such that for $i,j =1,2$ with $i\ne j$, if 
  $\{x_n: n \in \mathbb N \}$ and $\{y_n: n \in \mathbb N \}$ are two sequences in $X$ with $\{x_n: n \in 
  \mathbb N \}$ $\tau_i$-convergent to $p$ and $x_n \in g_j(n, y_n)$ for all $n \in \mathbb N$, then 
  $\{y_n: n \in \mathbb N \}$ is $\tau_i$-convergent to $p$.
  \end{itemize}
  \end{theorem}
  
  \begin{proof}
  $(1)\Rightarrow (2)$ follows directly from Theorem \ref{thm:pksemig1}, as $\{ p \} \cup \{x_n: n\in\mathbb 
  N\}$ is compact. 
  
  \medskip
  $(2) \Rightarrow (3)$. Let $(g_1, g_2)$ be a pairwise $g$-function satisfying (2). Let $\{x_n: n \in \mathbb 
  N \}$ and $\{y_n: n \in \mathbb N \}$ be two sequences in $X$ such that $\{x_n: n\in \mathbb N\}$ is 
  $\tau_i$-convergent to $p$ and $x_n \in g(n, y_n)$. Assume that $\{y_n: n \in \mathbb N \}$ is not 
  $\tau_i$-convergent $p$. Then $\{ y_n: n \in \mathbb N \}$ has a subsequence $\{y_{n_k}: k \in \mathbb N \}$ 
  such that $p \not \in \overline{\{ y_{n_k}: k \in \mathbb N \}}^{\tau_i}$. Put $H= \{ y_{n_k}: k \in \mathbb 
  N \}$. Since $\{x_{n_k}: k \in \mathbb N \}$ is $\tau_i$-convergent to $p$, then we can assume that $x_{n_k} 
  \not\in \overline{H}^{\tau_i}$ for all $k \ge 1$. Thus, by (2), there must be an $m \in \mathbb N$ such that 
  \[
  \left(\{ p \} \cup \{ x_{n_k}: k \in \mathbb N\} \right) \cap \left(\bigcup \{ g_j(m, x): x\in H\} \right) 
  = \emptyset.
  \]
  On the other hand, 
  \[
  x_{n_m} \in g_j(n_m, y_{n_m}) \subseteq g_j(m, y_{n_m}) \subseteq \bigcup \{ g_j(m, x): x\in H\}.
  \]
  A contradiction occurs.
  
  \medskip
  $(3) \Rightarrow (1)$. A proof has been given in \cite{LL}.
  \end{proof}
  
  In \cite{Gao:85}, $k$-semi-stratifiable spaces are defined in terms of $\sigma$-cushioned pair $k$-networks,
  which is different from (but equivalent to) that given in \cite{Lutzer:71}. Our next result, which just
  confirms that the same thing holds in the setting of bispaces, provides characterizations of a pairwise 
  $k$-semi-stratifiable bispace in terms of cushioned pair $k$-networks. 
  
  \begin{theorem} \label{thm:knetwork}
  A bispace $(X,\tau_1,\tau_2)$ is $\tau_1$-$k$-semi-stratifiable with respect to $\tau_2$ if, and only if, 
  it has a $\sigma$-$\tau_2$-cushioned pair $\tau_1$-$k$-network.
  \end{theorem}

  \begin{proof} 
  Necessity. Let $F_{12}: \mathbb N \times \tau_1 \to {\mathscr F}_2(X)$ be an operator 
  satisfying conditions (1)-(3) in Lemma \ref{lem:pksemi} such that $F_{12}$ is also monotone with respect to
  $n$. For each $n \in \mathbb N$, define $\Delta_n = \tau_1$ and
  \[
  {\mathscr B}_n = \left\{\left(B_U^1, B_U^2\right): B_U^1 = F_{12}(n, U), B_U^2=U \mbox{ and } U \in 
  \Delta_n \right\}.
  \]
  We claim that $\mathscr B_n$ is $\tau_2$-cushioned. Indeed, if $\mathscr U \subseteq \Delta_n$, 
  by condition (2) in Lemma \ref{lem:pksemi}, $F_{12}(n, U) \subseteq F_{12}\left(n, \bigcup 
  {\mathscr U}\right)$ for any $U \in \mathscr U$. It follows that
  \begin{eqnarray*}
  \overline{\bigcup \left\{B_U^1: U \in {\mathscr U} \right\}}^{\tau_2} &=&
  \overline{\bigcup \left\{F_{12}(n, U): U \in {\mathscr U} \right\}}^{\tau_2}\\
  &\subseteq& F_{12}\left(n, \bigcup {\mathscr U}\right)\\ 
  &\subseteq& \bigcup \{B_U^2: U \in \mathscr U\},
  \end{eqnarray*}
  which implies that each $\mathscr B_n$ is $\tau_2$-cushioned. Thus, $\mathscr B = \bigcup_{n\in \mathbb{N}}
  \mathscr B_n$ is $\sigma$-$\tau_2$-cushioned. Let $K$ be a $\tau_1$-compact subset of $X$ and $U\in\tau_1$
  with $K \subseteq U$. By condition (3) in Lemma \ref{lem:pksemi}, there must be some $m\in \mathbb N$ such 
  that $K \subseteq F_{12}(m, U) \subseteq U$. Then $\Delta' = \{U\} \subseteq \Delta_m$ and
  \[
  K \subseteq  \bigcup \{B_U^1: U \in \Delta'\} \subseteq \bigcup \{B_U^2: U \in \Delta'\} \subseteq U,
  \]
  which implies that $\mathscr B$ is also a pair $\tau_1$-$k$-network.
  
  \medskip
  Sufficiency. Let ${\mathscr B} = \bigcup_{n \in \mathbb N} {\mathscr B}_n$ be a $\sigma$-$\tau_2$-cushioned
  pair $\tau_1$-$k$-network, that is, $\mathscr B$ is a pair $\tau_1$-$k$-network and for each $n \in \mathbb 
  N$, ${\mathscr B}_n =\left\{ \left(B_\alpha^1, B_\alpha^2 \right): \alpha \in \Delta_n \right\}$ is a 
  $\tau_2$-cushioned family. Without loss of generality, for each $n \in \mathbb N$, we can assume that 
  ${\mathscr B}_n \subseteq {\mathscr B}_{n+1}$ and $\mathscr B_n$ is closed under finite union. Define 
  $F_{12}: \mathbb N \times \tau_1\to {\mathscr F}_2(X)$ such that for each $n\in \mathbb N$ and each $U \in 
  \tau_1$,
  \[
  F_{12}(n, U):= \overline{\bigcup \left\{ B_\alpha^1: B_\alpha^2 \subseteq U \mbox{ and } \alpha \in
  \Delta_n\right\}}^{\tau_2}.
  \]
  First of all, as $\mathscr B_n$ is $\tau_2$-cushioned, we have
  \[
  F_{12}(n, U) \subseteq \bigcup \left\{ B_\alpha^2: B_\alpha^2 \subseteq U \mbox{ and } \alpha \in
  \Delta_n\right\} \subseteq U.
  \]
  For each $x\in U$, as $\mathscr B$ is a pair $\tau_1$-$k$-network and $\{x\}$ is compact, there exist an 
  $m \in \mathbb N$ and an $\alpha \in \Delta_m$ such that $x \in B_\alpha^1 \subseteq B_\alpha^2 \subseteq 
  U$. It follows that $x\in F_{12}(m, U)$, which implies that $\bigcup_{n\in \mathbb N} F(n, U) =U$. It is 
  clear that if $U, V \in \tau_1$ with $U \subseteq V$, then $F_{12}(n, U) \subseteq F_{12}(n, V)$ for any 
  $n \in \mathbb N$. Finally, if $K$ is $\tau_1$-compact and $U\in \tau_1$ with $K \subseteq U$, similar 
  to what have done previously, there exist an $m \in \mathbb N$ and an $\alpha \in \Delta_m$ such that $K 
  \in B_\alpha^1 \subseteq B_\alpha^2 \subseteq U$. This implies that $K \subseteq F_{12}(m, U) \subseteq U$.
  Therefore, we have checked that $F_{12}$ is an operator satisfying conditions (1)-(3) in Lemma 
  \ref{lem:pksemi}.
  \end{proof}

  \begin{corollary}
  A bispace $(X,\tau_1,\tau_2)$ is pairwise $k$-semi-stratifiable if, and only if, it has a 
  $\sigma$-$\tau_j$-cushioned pair $\tau_i$-$k$-network for each pair of $i, j = 1, 2$ with $i\ne j$.
  \end{corollary}
  
  Let $(X, \tau_1, \tau_2)$ be a bispace and $x\in X$ be a point. A family ${\mathscr P}_x$ of subsets of $X$ 
  is called a \emph{$\tau_i$-cs-network at $x$} \cite{Guthrie:73} if for every sequence $\{x_n: n \in \N\}$ 
  that is $\tau_i$-convergent to $x$ and an arbitrary open neighborhood $U$ of $x$ in $(X, \tau_i)$, there 
  exist an $m\in\N$ and an element $P\in {\mathscr P}_x$ such that
  \[
  \{x\}\cup \{x_n: n\geqslant m\} \subseteq P \subseteq U.
  \]
  If each point $x$ in $X$ has a $\tau_i$-cs-network ${\mathscr P}_x$, then $\bigcup_{x\in X}{\mathscr P}_x$ 
  is called \emph{$\tau_i$-cs-network for $(X, \tau_i)$}.
  
  \medskip
  In \cite{Gao:85}, Gao characterized $k$-semi-stratifiable spaces in terms of $cs$-networks. At the end of 
  this section, we establish a similar result in the setting of bispaces.
  
  \begin{theorem} \label{thm:csnetwork}
  Let $(X, \tau_1, \tau_2)$ be a bispace such that $(X, \tau_i)$ is Hausdorff for $i=1,2$. Then $(X, \tau_1, 
  \tau_2)$ is pairwise $k$-semi-stratifiable if, and only if, for any $i, j=1,2$ with $i\ne j$, there is 
  an operator $F_{ij}: \mathbb{N}\times \tau_i \rightarrow {\mathscr F}_j(X)$ satisfying
  \begin{enumerate}
  \item[{\rm (1)}] $U=\bigcup_{n\in \mathbb{N}}F_{ij}(n,U)$ for all $U \in \tau_i$;
  	
  \item[{\rm (2)}] if $U, V\in \tau_i$ with $U \subseteq V$, then $F_{ij}(n,U)\subseteq F_{ij}(n,V)$ for all 
  $n\in \mathbb{N}$;
  	
  \item[{\rm (3)}] for each $U\in \tau_i$, $\{ F_{ij}(n,U): n \in \mathbb N\}$ is a $\tau_i$-$cs$-network at 
  every point of $U$.
  \end{enumerate}
  In addition, the operator $F_{ij}$ can be required to be monotone with respect to $n$, that is, $F_{ij}(n,U) 
  \subseteq F_{ij}(n+1,U)$ for all $n\in \mathbb{N}$ and all $U\in \tau_i$.  
  \end{theorem}
  
  \begin{proof}
  The necessity is trivial by Lemma \ref{lem:pksemi}, as $\{x \} \cup \{x_n: n \in \mathbb N\}$ is 
  $\tau_i$-compact for any sequence $\{x_n: n \in \mathbb N\}$ in $X$ which is $\tau_i$-convergent to a 
  point $x \in X$.
  
  \medskip
  Sufficiency. Suppose that for any $i, j=1,2$ with $i\ne j$, there is an operator $F_{ij}: \mathbb{N}\times 
  \tau_i \rightarrow {\mathscr F}_j(X)$ satisfying conditions (1)--(3) above and monotonicity. We only need
  to verify condition (3) in Lemma \ref{lem:pksemi}. First, note that these conditions imply that each point 
  $x$ is a $G_\delta$-set in both $(X, \tau_1)$ and 
  $(X, \tau_2)$. Suppose that there are a $\tau_i$-compact set $K$ and a $U \in \tau_i$ with $K \subseteq 
  U$, but $K \not \subseteq F_{ij}(n, U)$ for any $n \in \mathbb N$. Then, there is a sequence $\{x_n: n 
  \in \mathbb N\} \subseteq K$ such that $x_n \not \in F_{ij}(n, U)$ for any $n \in \mathbb N$. Since $K$
  is $\tau_i$-compact and points are $G_\delta$ in $(X, \tau_i)$, then $\{x_n: n \in \mathbb N\}$ must 
  have a subsequence $\{x_{n_k}: k \in \mathbb N\}$ which is $\tau_i$-convergent to a point $x\in K$ and 
  $x_{n_k} \not \in F_{ij}(n_k, U)$ for all $k \in \mathbb N$. By condition (3) above, there exist an 
  $m_0\in \mathbb N$ and an $n_0 \in \mathbb N$ such that
  \[
  \{x\}\cup \{x_{n_k}: k\geqslant m_0\} \subseteq F_{ij}(n_0, U) \subseteq U.
  \]
  It follows that for any $k\in \mathbb N$ with $k \ge m_0$ and $n_k \ge n_0$, we have
  \[
  x_{n_k} \in F_{ij}(n_0, U) \subseteq F_{ij}(n_k, U).
  \]
  Apparently, this contradicts with the choice of $x_{n_k}$. We have verified that condition (3) in
  Lemma \ref{lem:pksemi} is satisfied, and thus  $(X, \tau_1, \tau_2)$ is  pairwise $k$-semi-stratifiable.
  \end{proof}

  Note that Theorems \ref{thm:knetwork} and \ref{thm:csnetwork} may shed some light on relationships 
  between pairwise $k$-semi-stratifiability and the other generalized metric properties of bispaces studied
  in \cite{Mu} and other places.

  \section{When is a pairwise $k$-semi-stratifiable bispaces\\ pairwise
  straitifiable?} \label{sec:relation}
  
  In this section, we consider the problem when a pairwise $k$-semi-stratifiable bispace is pairwise 
  stratifiable. An open question posed in \cite{LL} is completely solved, and some results in 
  \cite{Cao:99}, \cite{Cao:05} and \cite{Gao:85} are extended to the setting of bispaces.
 
  \medskip
  Recall that a topological space $(X, \tau)$ is said to be \emph{Fr\'{e}chet}, if for every nonempty subset
  $A \subseteq X$ and every point $x\in \overline{A}^{\tau}$, there is a sequence $\{x_n: n \in \mathbb 
  N\} \subseteq A$ such that  $\{x_n: n \in \mathbb N\}$ converges to $x$. 
  
  \medskip
  In a recent paper \cite{LL}, Li 
  and Lin posed the following open question (see \cite[Question 3.4]{LL}).
  
  \begin{question}[\cite{LL}] \label{ques:lilin}
  Is a pairwise $k$-semi-stratifiable bispace $(X, \tau_1, \tau_2)$ pairwise stratifiable if 
  $(X, \tau_i)$ is a Fr\'{e}chet space for each $i = 1, 2$?
  \end{question}
  
  Our next theorem answers Question \ref{ques:lilin} affirmatively. Note that our theorem also extends a
  result in \cite{Gao:85} to the setting of bispaces.
  
  \begin{theorem}
  Let $(X, \tau_1, \tau_2)$ be a pairwise $k$-semi-stratifiable bispace. If both $(X, \tau_1)$ and 
  $(X, \tau_2)$ are Fr\'{e}chet spaces, then $(X, \tau_1, \tau_2)$ is pairwise stratifiable.
  \end{theorem}
  
  \begin{proof}
  In the light of Corollary \ref{coro:bidecom}, we need to show that $(X, \tau_1, \tau_2)$ is pairwise
  monotonically normal. For any fixed $i, j=1,2$ with $i\ne j$, let $F_{ij}: \mathbb{N}\times \tau_i 
  \rightarrow {\mathscr F}_j(X)$ be an operator that is monotone with respect to $n$ and satisfies (1)-(3) 
  in Lemma \ref{lem:pksemi}. For each pair $(H, K)$ of disjoint subsets of $X$ such that $H \in {\mathscr 
  F}_i(X)$ and $K \in {\mathscr F}_j(X)$, define $D_{ij}(H,K)$ by
  \[
  D_{ij}(H,K) := {\rm int}_{\tau_j}\left(\bigcup_{n\in \mathbb N} \left(F_{ji}(n, X \smallsetminus K) 
  \smallsetminus F_{ij}(n, X \smallsetminus H)\right)\right)
  \]
  Next, we shall verify that $D_{ij}(\cdot, \cdot)$ satisfies all conditions in the definition of a pairwise
  monotonically normal bispace.
  
  \medskip
  Clearly, $D_{ij}(H, K) \in \tau_j$ and $D_{ij}(\cdot, \cdot)$ satisfies condition (ii) in the definition 
  of a pairwise monotonically normal bispace. Also note that $D_{ij}(H,K) \subseteq X 
  \smallsetminus K$ holds trivially.
	
  \medskip
  \noindent
  \textbf{Claim 1}. $H \subseteq D_{ij}(H,K)$. 
  \begin{proof}[Proof of Claim 1]
  Suppose that there is a point $x_0 \in H \smallsetminus D_{ij}(H,K)$. Then 
  \begin{eqnarray*}
  x_0 &\in& X\smallsetminus {\rm int}_{\tau_j} \left(\bigcup_{n\in\mathbb N} \left(F_{ji}(n,X\smallsetminus K) 
  \setminus F_{ij}(n, X \smallsetminus H)\right)\right)\\ 
     &=& 
  \overline{X\smallsetminus\bigcup_{n\in\mathbb N} \left(F_{ji}(n, X \smallsetminus K) \smallsetminus F_{ij}
  (n, X \smallsetminus H)\right)}^{\tau_j}.
  \end{eqnarray*}
  Since $(X, \tau_j)$ is a Fr\'{e}chet space, there is a sequence $\{x_n: n \in \mathbb N\}$ such that
  \[
  \{x_n: n \in \mathbb N\} \subseteq X \smallsetminus \bigcup_{n\in \mathbb N} \left(F_{ji}(n, X 
  \smallsetminus K) \smallsetminus F_{ij}(n, X \smallsetminus H)\right)
  \] 
  and $\{x_n: n \in \mathbb N\}$ is $\tau_j$-convergent to $x_0$. Note that $x_0 \in H$ implies that $x_0 
  \in X \smallsetminus K$. Thus, there is an $m\in \mathbb N$ such that 
  \[
  C = \{x_0 \} \cup \{x_n: n \ge m\} \subseteq X \smallsetminus K. 
  \]
  By (3) in Lemma \ref{lem:pksemi}, there is an $m' \ge m$ such that $C \subseteq F_{ji}(m', X 
  \smallsetminus K)$. On the other hand, note that there must be some $p \ge m'$ such that $x_p \not \in 
  F_{ij}(m', X \smallsetminus H)$. Otherwise, as $F_{ij}(m', X \smallsetminus H)$ is $\tau_j$-closed, we 
  conclude that
  \[
  x_0 \in F_{ij}(m', X \smallsetminus H) \subseteq X \smallsetminus H,
  \] 
  which contradicts with the fact $x_0\in H$. It follows that
  \begin{eqnarray*}
  x_p &\in& F_{ji}(m', X \smallsetminus K) \smallsetminus F_{ij}(m', X \smallsetminus H)\\ 
  &\subseteq&  \bigcup_{n\in \mathbb N} \left(F_{ji}(n, X \smallsetminus K) \smallsetminus F_{ij}
  (n, X \smallsetminus H)\right).
  \end{eqnarray*}
  This certainly contradicts with the selection of $\{x_n: n \in \mathbb N\}$. Hence, Claim 1 has been 
  verified.
  \end{proof}
	
  \medskip
  \noindent
  \textbf{Claim 2}. $\overline{D_{ij}(H,K)}^{\tau_i} \subseteq X \smallsetminus K$.
  \begin{proof}[Proof of Claim 2]	
  Suppose that there is a point $x_0 \in \overline{D_{ij}(H,K)}^{\tau_i} \cap K$. Since $(X,\tau_i)$ is a 
  Fr\'{e}chet space, there is a sequence $\{x_n: n \in \mathbb N\} \subseteq D_{ij}(H,K)$ such that $\{x_n: 
  n\in \mathbb N\}$ is $\tau_i$-convergent to $x_0$. Note that $x_0 \in K$ implies $x_0\in X \smallsetminus 
  H$. Thus, there exists an $m \in \mathbb N$ such that 
  \[
  \{x_0\} \cup \{x_n: n \ge m\} \subseteq X \smallsetminus H.
  \] 
  By condition (3) in Lemma \ref{lem:pksemi}, there exists some $p \ge m$ such that 
  \[
  \{x_0\} \cup \{x_n: n \ge m\} \subseteq F_{ij}(p, X \smallsetminus H).
  \] 
  By the selection of $\{x_n: n \in \mathbb N\}$, we know that
  \[
  \{x_n: n \ge m\} \subseteq \bigcup_{n\in \mathbb N} \left(F_{ji}(n, X \smallsetminus K) \smallsetminus 
  F_{ij}(n, X \smallsetminus H)\right),
  \]
  which implies that
  \[
  \{x_n: n \ge m\} \subseteq \bigcup_{n= 1}^{p-1} \left(F_{ji}(n, X \smallsetminus K) \smallsetminus 
  F_{ij}(n, X \smallsetminus H)\right).
  \]
  It follows that there are a $q \in \mathbb N$ with $1\le q <p$ and a subsequence $\{x_{n_k}: k \ge m\}$ 
  of $\{x_n: n \in \mathbb N\}$ such that 
  \[
  \{x_{n_k}: k \ge m\} \subseteq F_{ji}(q, X \smallsetminus K) 
  \smallsetminus F_{ij}(q, X \smallsetminus H).
  \] 
  We conclude that $x_0 \in F_{ji}(q, X \smallsetminus K) 
  \subseteq X \smallsetminus K$. This contradicts with the selection of $x_0$, and thus Claim 2 has been 
  verified.
  \end{proof}
	
  Combining Claims 1 and 2, we see that $D_{ij}(\cdot, \cdot)$ also satisfies condition (i) in the 
  definition of a pairwise monotonically normal bispace.
  \end{proof}
  
  Let $(X,\tau)$ be a topological space, and let $x \in X$ be a point. The collection of neighborhoods of $x$ 
  in $(X,\tau)$ is denoted by ${\mathscr N}(\tau, x)$. We shall consider the following \emph{${\mathscr 
  G}(\tau, x)$-game} played in $(X,\tau)$ between two players: $\alpha$ and $\beta$. Player $\alpha$ goes first 
  and chooses a point $x_1 \in X$. Player $\beta$ then responds by choosing $U_1 \in {\mathscr N}(\tau, x)$. 
  Following this, $\alpha$ must select another (possibly the same) point $x_2 \in U_1$ and in turn $\beta$ 
  must again respond to this by choosing (possibly the same) $U_2 \in {\mathscr N}(\tau, x)$. The players 
  repeat this procedure infinitely many times, and produce a play $(x_1, U_2, x_2, U_2, \cdots, x_n, U_n,
  \cdots)$ in the ${\mathscr G}(\tau, x)$-game, satisfying $x_{n+1} \in U_n$ for all $n \in \mathbb N$. We
  shall say that $\beta$ \emph{wins} a play $(x_1, U_2, x_2, U_2,\cdots,x_n,U_n,\cdots)$ if the sequence 
  $\{ x_n: n \in \mathbb N \}$ has a cluster point in $X$. Otherwise, $\alpha$ is said to \emph{have won} 
  the play. By a \emph{strategy} $s$ for $\beta$, we mean a `rule' that specifies each move of $\beta$ in 
  every possible situation. More precisely, a strategy $s$ for $\beta$ is an ${\mathscr N}(\tau, x)$-valued 
  function. We shall call a finite sequence $\{ x_1, x_2,..., x_n \}$ or an infinite sequence $\{ x_1, x_2,...\}$ 
  an \emph{$s$-sequence} if $x_{i+1} \in s(x_1, x_2,...,x_i)$ for each $i$ such that $1 \leq i < n$ or $x_{n+1} 
  \in s(x_1, x_2,...,x_n)$ for each $n \in \mathbb N$. A strategy $s$ for player $\beta$ is called \emph{a 
  winning strategy} if each infinite $s$-sequence has a cluster point in $X$. Finally, we call $x$ a 
  \emph{$\mathscr G(\tau)$-point} if player $\beta$ has a winning strategy for the ${\mathscr G}(\tau,x)$-game. 
  In addition, if every point of $X$ is a $\mathscr G(\tau)$-point, then $(X,\tau)$ is called a 
  \emph{$\mathscr G$-space} \cite{Bouziad:93}. The notion of $\mathscr G$-spaces is a common generalization 
  of the concepts of $q$-spaces in \cite{Michael:64} and $W$-spaces in \cite{Gr:76}. 
  
  \medskip
  The next result extends \cite[Theorem 2.2.8]{Cao:99} and \cite[Theorem 3.2]{Cao:05} to the setting of bispaces.
  
  \begin{theorem}
  Let $(X, \tau_1, \tau_2)$ be a pairwise $k$-semi-stratifiable bispace. If both $(X, \tau_1)$ and $(X, 
  \tau_2)$ are regular and $\mathscr G$-spaces, then $(X, \tau_1, \tau_2)$ is pairwise stratifiable.
  \end{theorem}
  
  \begin{proof}
  Let $(g_1, g_2)$ be a pairwise $g$-function as described in condition (3) of Theorem \ref{thm:pksemig2}. For
  $i, j=1,2$ with $i\ne j$, we define $G_{ij}: \mathbb N \times {\mathscr F}_i(X) \to \tau_j$ such that for
  each $n \in \mathbb N$ and each $H \in \mathscr F_i(X)$,
  \[
  G_{ij}(n, H) = \bigcup \left\{ g_j(n, x) : x \in H\right\}.
  \]
  Clearly, if $H, K \in {\mathscr F}_i(X)$ with $H\subseteq K$, then $G_{ij}(n,H)\subseteq G_{ij}(n,K)$ for 
  all $n\in \mathbb N$. Furthermore, it is easy to see that $H\subseteq \bigcap_{n\in \mathbb N} 
  \overline{G_{ij} (n,H)}^{\tau_i}$ for all $H\in {\mathscr F}_i(X)$.
  	
  \medskip
  Suppose that there are a point $x \in X$ and a $\tau_i$-closed subset $H$ in $X$ with $x \not \in H$, but 
  $x \in \overline{G_{ij}(n,H)}^{\tau_i}$ for every $n \in \mathbb N$. First, we choose some $\tau_i$-open 
  neighborhood $U$ of $x$ such that $\overline{U}^{\tau_i} \cap H = \emptyset$. Since $(X,\tau_i)$ is a 
  $\mathscr G$-space, $\beta$ has a winning strategy $s$ for the ${\mathscr G}(\tau_i, x)$-game. 
  Let $\alpha$'s first move be $x_1$. By our assumption and the definition of $G_{ij}(\cdot, \cdot)$, there 
  must exist some point $y_1 \in H$ such that $s(x_1) \cap U \cap g_j(1, y_1) \ne \emptyset$.
  Inductively, we can obtain two sequences $\{x_n: n \in \mathbb N\}$ and $\{y_n: n \in \mathbb N\}$ in $X$ 
  such that for each $n \in \mathbb N$, $y_n \in H$ and
  \[
  x_{n+1} \in U \cap g_j(n+1, y_{n+1}) \cap \left(\bigcap_{1 \leq j \leq n, \atop 1\leq i_1 \leq...\leq i_j
  \leq n} s(x_{i_1},...,x_{i_j}) \right).
  \]
  It follows that each subsequence of $\{ x_n: n \in \mathbb N\}$ is an $s$-sequence in $(X,\tau_i)$, and 
  thus has an cluster point in $(X, \tau_i)$. Since each point of $X$ is a $G_\delta$-point in $(X,\tau_i)$, 
  then $\{ x_n: n \in \mathbb N\}$ must have a convergent subsequence, saying $\{x_{n_k}: k\in \mathbb N\}$, 
  in $(X, \tau_i)$. Suppose that $\{ x_{n_k}: k \in \mathbb N\}$ is $\tau_i$-convergent to some point $ x_* \in
  \overline{U}^{\tau_i}$. Then, by condition (3) in Theorem \ref{thm:pksemig2}, and the construction of $\{x_n: 
  n \in \mathbb N\}$ and $\{y_n: n \in \mathbb N\}$ in the above, $\{ y_{n_k}: k\in \omega\}$ is also 
  $\tau_i$-convergent to $x_*$, and thus $x_* \in H$. It follows that $x_* \in \overline{U}^{\tau_i} \cap H$. 
  We have derived a contradiction. Therefore, $x \not \in \overline{G_{ij}(n,H)}^{\tau_i}$ for some $n \in 
  \mathbb N$. We have verified that $H =\bigcap_{n\in \mathbb N} \overline{G_{ij}(n,H)}^{\tau_i}$ for all $H
  \in {\mathscr F}_i(X)$ and thus $(X, \tau_1, \tau_2)$ is pairwise stratifiable.
  \end{proof}

  \section{Quasi-pseudo-metrizability of topological ordered spaces} \label{sec:quasi}
  
  In \cite{KM}, K\"{u}nzi and Mushaandja posed the following open problem (refer to \cite[Problem 1]{KM}).
  
  \begin{problem}[\cite{KM}] \label{prob:KM}
  If $(X, \tau, \preccurlyeq)$ is a topological ordered $C$- and $I$-space such that the topology $\tau$ is 
  metrizable, is the associated bitopological space $(X,\tau^{\flat},\tau^{\natural})$ quasi-pseudo-metrizable?
  \end{problem}
  
  It was shown that if the topology $\tau$ is separable metrizable, then $(X,\tau^{\flat},\tau^{\natural})$ 
  is quasi-pseudo-metrizable. This result gives a partial affirmative answer to Problem \ref{prob:KM} in the 
  class of separable metrizable topological ordered $C$- and $I$-spaces. In this section, we provide another 
  partial affirmative answer to this problem in the class of ball transitive and metrizable topological 
  ordered $C$- and $I$-spaces. To this purpose, we first introduce the concept of ball transitivity.
  
  \medskip
  An ordered metric space $(X, \rho, \preccurlyeq)$ is said to be \emph{ball transitive} \cite{Ri} if there 
  exists an $n\in \mathbb N$ such that whenever $x\preccurlyeq y$, then $B_\rho(x, \frac{\epsilon}{n}) 
  \subseteq d(B_\rho(y,\epsilon))$ and $B_\rho(y, \frac{\epsilon}{n}) \subseteq i(B_\rho(x,\epsilon))$ hold 
  for any $\epsilon>0$. We call a metrizable topological ordered space $(X, \tau, \preccurlyeq)$ \emph{ball 
  transitive} provided that there is a metric $\rho$ compatible with $\tau$ such that $(X,\rho,\preccurlyeq)$ 
  is ball transitive. 
  
  \begin{remark}
  (i) Let $X$ be the space of continuous real-valued functions on the interval $[0, 1]$. Let $\preccurlyeq$ be 
  the pointwise order on $X$ and $\rho$ be the metric defined by the sup-norm. It is well known that $(X,\rho)$
  is not separable, but it was shown in \cite{Ri} that is ball transitive with $n=1$.
  	
  \medskip
  (ii) Let $A$ be the open first quadrant of $\mathbb R^2$, i.e., 
  \[
  A=\{(x,y) \in \mathbb R^2: x>0,y>0 \}. 
  \]
  Consider the following subset $X$ of $\mathbb R^2$,
  \[
  X = \left([-1, 1] \times [-1,1] \right) \setminus A,
  \]
  equipped with the Euclidean metric $\rho$ and the pointwise order $\preccurlyeq$. It is clear that $(X, 
  \rho, \preccurlyeq)$ is separable. However, it was shown in \cite{Ri} that $(X, \rho, \preccurlyeq)$ is not 
  ball transitive.
  \end{remark}
  
  Recall that a topological space $(X,\tau)$ is said to be a \emph{$\gamma$-space} \cite{Gr}, provided that 
  there is a $g$-function $g: \mathbb N \times X \to \tau$ such that if $x_n \in g(n,x)$ and $g_n \in 
  g(n,x_n)$ for all $n \in \mathbb N$ then $x$ is a cluster point of the sequence $\{y_n: n \in \mathbb N\}$.
  Herein, we call such a function $g$ a \emph{$\gamma$-function} for $(X,\tau)$.
  
  \begin{theorem} \label{thm:ball}
  Let $(X,\tau,\preccurlyeq)$ be a metrizable topological ordered $I$-space. If $(X,\tau,\preccurlyeq)$ is 
  ball transitive, then $(X,\tau^{\flat})$ and $ (X,\tau^{\natural})$ are $\gamma$-spaces.
  \end{theorem}

  \begin{proof} Let $\rho$ be a metric compatible with $\tau$ such that $(X, \rho,\preccurlyeq)$ is ball 
  transitive. Then, there is a $k \in \mathbb N$ such that whenever $x\preccurlyeq y$, $B_\rho(x, 
  \frac{\epsilon}{k}) \subseteq d(B_\rho(y,\epsilon))$ and $B_\rho(y, \frac{\epsilon}{k}) \subseteq 
  i(B_\rho(x,\epsilon))$ hold for any $\epsilon>0$. For each $x\in X$, let ${\mathscr B}(x) = \left\{ B_\rho(x,\frac{1}{2^{n}}):n\in \mathbb N\right\}$. Since $(X,\tau, \preccurlyeq)$ is a metrizable $I$-space, 
  then for any $x\in X$, ${\mathscr U}(x) =\left\{ d \left(B_\rho(x,\frac{1}{2^{n}}) \right): n \in \mathbb 
  N\right\}$ is a countable base of open neighbourhoods for $\tau^{\flat}$ at $x$ and ${\mathscr V}(x) = 
  \left\{ i\left(B_\rho(x,\frac{1}{2^{n}}) \right): n \in \mathbb N \right\}$ is a countable base of open 
  neighbourhoods for $\tau^{\natural}$ at $x$, respectively.
  
  \medskip
  Next, define a $g$-function $g_1: {\mathbb N}\times X\rightarrow \tau^{\flat}$ by letting $g_1(n,x) =
  d\left(B_\rho(x,\frac{1}{k4^{n}})\right)$ for each $x\in X$ and $n\in \mathbb N$. We verify that
  $g_1: {\mathbb N}\times X\rightarrow \tau^{\flat} $ is a $\gamma$-function for $(X,\tau^{\flat})$. 
  To this end, let $\{x_n: n \in \mathbb N\}$ and $\{y_n: n \in \mathbb N\}$ be two sequences in $X$ such 
  that $x_n\in g_1(n,x)$ and $y_n\in g_1(n,x_n)$ for all $n \in \mathbb N$. Without loss of generality, 
  we may require  $y_n\neq x$ for any $n\in \mathbb N$. Suppose that $x$ is not a 
  $\tau^\flat$-cluster point of $\{y_n: n\in \mathbb N\}$. Then there exists an $m\in \mathbb N$ such that
  \[
  d\left(B_\rho(x,\frac{1}{2^{m}})\right) \subseteq X \setminus \overline{\{ y_n: n\in \mathbb N\}}^{\tau^{\flat}}. 
  \] 
  For each $n\geqslant m$, as $x_n\in g_1(n,x)=d(B_\rho(x,\frac{1}{k4^{n}}))$, there exists a $t_n\in B_\rho
  (x,\frac{1}{k4^{n}})$ such that $x_n\preccurlyeq t_n$. By the ball transitivity of $(X,\rho,\preccurlyeq)$, 
  we have $B_\rho(x_n, \frac{1}{k4^m})\subseteq d\left(B_\rho(t_n, \frac{1}{4^m})\right)$, which implies that
  \[
  B_\rho(x_n, \frac{1}{k4^m})\subseteq d\left(B_\rho(t_n, \frac{1}{4^m})\right) \subseteq 
  d\left(B_\rho(x, \frac{1}{2^m})\right)
  \]
  for all $n \ge m$. It follows that $y_n \in d\left(B_\rho(x, \frac{1}{2^m})\right)$ for all $n \ge m$.
  This is contradiction. Hence, $x$ is not a $\tau^\flat$-cluster point of $\{y_n: n\in \mathbb N\}$, which 
  verifies that $g_1$ is a $\gamma$-function for $(X,\tau^{\flat})$.
  
  \medskip
  Finally, define a $g$-function $g_2: {\mathbb N}\times X\rightarrow \tau^{\natural}$ by letting $g_2(n,x)= i\left(B_\rho(x,\frac{1}{k4^{n}})\right)$ for each $x\in X$ and $n\in \mathbb N$. In the way similar to 
  what we have done previously, we can prove that $g_2$ is a $\gamma$-function for $(X,\tau^{\natural})$. 
  Therefore, both $(X,\tau^{\flat})$ and $ (X,\tau^{\natural})$ are $\gamma$-spaces.
  \end{proof}

  \begin{lemma}[\cite{Marin:09}] \label{lem 2.1} 
  A bispace $(X,\tau_1, \tau_2)$ is quasi-pseudo-metrizable if, and only if, $(X,\tau_1, \tau_2)$ is pairwise 
  stratifiable and $(X,\tau_i)$ is a $\gamma$-space for $i=1, 2$.
  \end{lemma}
  
  The following result provides a partial answer to Problem \ref{prob:KM}.
  
  \begin{theorem} \label{thm:qpmetric}
  Let $(X,\tau,\preccurlyeq)$ be a topological ordered $C$- and $I$-space such that the topology $\tau$ is 
  metrizable. If $(X,\tau,\preccurlyeq)$ is ball transitive, then $(X,\tau^{\flat},\tau^{\natural})$ is 
  quasi-pseudo-metrizable.
  \end{theorem}

  \begin{proof} 
  First, by Theorem \ref{thm:ball}, both $(X,\tau^{\flat})$ and $(X, \tau^{\natural})$ are $\gamma$-spaces. 
  Furthermore, by \cite[Theorem 1]{KM}, $(X,\tau^{\flat},\tau^{\natural})$ is a pairwise stratifiable 
  bispace. Hence, it follows from \cite[Theorem 4]{Marin:09} that $(X,\tau^{\flat},\tau^{\natural})$ is 
  quasi-pseudo-metrizable.
  \end{proof}
  
  Let $(X,\tau,\preccurlyeq)$ be a topological ordered $C$-space. It was shown in \cite{KM} that if $\tau$ 
  is a stratifiable topology, then $(X,\tau^{\flat},\tau^{\natural})$ is pairwise stratifiable. In addition, 
  it was shown in \cite{Kedian-Fucai:13} that if $\tau$ is a semi-stratifiable (resp. monotonically normal) 
  topology, then $(X,\tau^{\flat},\tau^{\natural})$ is pairwise semi-stratifiable (resp. monotonically normal).
  In the light of these results, we conclude this paper by posing the following open question.
  
  \begin{question}
  Let $(X,\tau,\preccurlyeq)$ be a topological ordered $C$-space. If $\tau$ is a $k$-semi-stratifiable topology,
  must $(X,\tau^{\flat},\tau^{\natural})$ be pairwise $k$-semi-stratifiable?
  \end{question}
  
  \medskip
  \centerline{\sc Acknowledgement}
  
  \bigskip
  The authors acknowledge the privilege of having seen \cite{LL} before its publication. 
  
  
\end{document}